\documentclass{naco}
  \usepackage{paralist}
  \usepackage{graphics} 
\usepackage{epstopdf}

\usepackage{mathrsfs}
\usepackage{amssymb}
\usepackage{amsmath}
\usepackage{cite}
\usepackage{graphicx}
\usepackage{subfig}
\usepackage{color}
\usepackage{float}
\usepackage{multicol}
\usepackage{multirow}
\usepackage{soul}
\usepackage{bm}
\usepackage{enumerate}

\def\eq#1{\begin{equation}#1\end{equation}}

 \def\qed{ \rule{.1in}{.1in}}

\def\currentvolume{x}
 \def\currentissue{x}
  \def\currentyear{201x}
   \def\currentmonth{X}
    \def\ppages{000--000}
     \def\DOI{X}

\newtheorem{theorem}{Theorem}
\newtheorem{definition}{Definition}
\newtheorem{lemma}{Lemma}

\newcommand{\diag}{{\rm diag\;}}
\newcommand{\col}{{\rm col\;}}

\title[Resilient Convex Combination for consensus-based algorithms] 
      {A Resilient Convex Combination for consensus-based distributed algorithms}

\author[X. Wang, S. Mou, S. Sundaram]{}



 \email{wang3156@purdue.edu}
 \email{mous@purdue.edu}
 \email{sundara2@purdue.edu}

\thanks{$^*$ Corresponding author: Shaoshuai Mou}

\begin{document}
\maketitle

\centerline{\scshape Xuan Wang and Shaoshuai Mou$^*$}
\medskip
{\footnotesize
 \centerline{School of Aeronautics and Astronautics }
   \centerline{Purdue University, West Lafayette, IN, USA}
} 

\medskip

\centerline{\scshape Shreyas Sundaram}
\medskip
{\footnotesize
 \centerline{School of Electrical and Computer Engineering}
   \centerline{Purdue University, West Lafayette, IN, USA}

}

\bigskip


\begin{abstract}
Consider a set of vectors in $\mathbb{R}^n$, partitioned into two classes: normal vectors and malicious vectors. The number of malicious vectors is bounded but their identities are unknown.
The paper provides a way for achieving a resilient convex combination, which is a convex combination of only normal vectors. Compared with existing approaches based on Tverberg points, the proposed method based on the intersection of convex hulls has lower computational complexity. Simulations suggest that the proposed method can be applied to resilience for consensus-based distributed algorithms against Byzantine attacks.
\end{abstract}

\section{INTRODUCTION}
Besides his significant research contribution to optimal control\cite{AP71ASME}, adaptive systems \cite{AB86MIT} and communications \cite{ZRB91TOC}, Professor Brian David Outram Anderson has recently focused on developing distributed algorithms for multi-agent networks \cite{MAB08TAC,BSUA16NACO,JAB07SIAM} based on the idea of consensus. Consensus-based distributed algorithms enable all agents in a network to reach an agreement regarding a certain quantity of interest, which could be an unconstrained value \cite{AJA03TAC}, a solution to a group of linear equations \cite{SJA15TAC}, or a constant for optimizing an objective function \cite{AA09TAC}. Success of these updates heavily depends on the utilization of convex combinations of nearby neighbors' states. When one or more agents become malicious under cyber-attacks, false information will be injected into the convex combination and usually lead to failures of consensus-based distributed algorithms \cite{Shreyas11TAC,Shreyas13Selected, Shreyas11TAC16Arxiv,AS18arxiv}. Considering the fact that many multi-agent networks in practice such as distributed power grids or robotic networks, are large-scale and often operate in open and hostile environments, the exposure to cyber-attacks is  inevitable \cite{FFF13TAC}. Moreover, in a fully distributed scenario, the lack of global information makes it almost impossible to identify or isolate those malicious agents, especially when the cyber-attack is very sophisticated such as Byzantine attack \cite{LRM82ACM}. Although significant progress has recently been achieved by a combination of cyber and system-theoretic approaches in \cite{ZMS10TCS,FAF12TAC, NLG12ACM}, these methods are either computationally expensive, assume the network topology to be fully connected, or require the normal nodes to be aware of nonlocal information such as independent paths between themselves and other nodes. Recognition of this has motivated us to achieve a \emph{resilient convex combination}, which refers to the convex combination of normal states that have not been manipulated by cyber-attacks, only knowing the upper bound to the number of malicious agents. Such a resilient convex combination could be achieved by recently developed methodologies based on Tverberg points \cite{mendes2015DC,LN14ACM,N14ICDCN,HS17TOR,HS16ICRA}, which are however computationally expensive. Thus one major goal of this paper is to develop an algorithm with low computational complexity for achieving a resilient convex combination. We will also apply the resilient convex combination in providing safety for consensus-based distributed algorithms in the adversarial environment.

The rest of this paper is organized as follows. In Section \ref{ProbF}, we formulate the problem of interest, and present a method for achieving a resilient convex combination in Section \ref{Sec_RCC}. The method is based on intersection of convex hulls and can be implemented by solving an optimization problem with low-computational complexity.
In Section \ref{Sec_App}, we apply the proposed resilient convex combination to consensus-based distributed algorithms and provide
numerical simulations to validate its effectiveness.  Finally, we conclude the paper in Section \ref{Conc}.

\smallskip
\noindent{ \emph{Notation}:}
Throughout this paper, we let ${\bf 1}_r$ denote a vector in $\mathbb{R}^r$ with all its components equal to 1; let $I_r$ denote the $r\times r$ identity matrix.
The transpose and kernel of a matrix $M$ are denoted by $M'$ and $\ker M$, respectively. For a square matrix $M$, by $M> 0$ and $M\geq 0$, we mean that the matrix $M$ is positive definite and positive semi-definite, respectively. For a vector $\beta$, by $\beta> 0$ and $\beta\geq 0$, we mean that each entry of vector $\beta$ is positive and non-negative, respectively.
We let $\diag\{A_1,A_2,\cdots,A_r\}$ denote the block diagonal matrix with $A$ the $i$th diagonal block entry, $i=1,2,\cdots,r$. Let $\otimes$ denote the Kronecker product and $\|\cdot\|_2$ denote the $2$-norm. 

\section{Problem Formulation} \label{ProbF}

Let $x_{\mathcal{A}}=\{x_1,x_2,\cdots,x_m\}$ denote a set of vectors in $\mathbb{R}^n$, where ${\mathcal{A}}=\{1,2,...,m\}$. Suppose one knows that at most a number of $\kappa$ vectors in $x_{\mathcal{A}}$ are malicious, but the labels of malicious vectors are not known. Then there are at least a number of $p=m-\kappa$ normal vectors in $x_{\mathcal{A}}$.  Suppose one knows a subset $\bar{\mathcal{A}}\subset \mathcal{A}$, which is empty or only contains labels of normal vectors in $x_{\mathcal{A}}$ but
\eq{ |\bar{\mathcal{A}}|=\sigma \leq p.}
The \textbf{problem of interest} is to develop an  algorithm with low-computational complexity to achieve a \emph{resilient convex combination}, which is defined as follows
\begin{definition} (resilient convex combination)
A vector is a resilient convex combination of $x_{\mathcal{A}}$, if it is a convex combination of at least $p$ normal vectors in $x_{\mathcal{A}}$, where $p$ is an known lower bound of normal vectors in  $x_{\mathcal{A}}$.
\end{definition}

The problem is trivial when $\kappa=0$ or $\sigma \geq p$, for which a resilient convex combination simply becomes a convex combination of $x_{\mathcal{A}}$ and $x_{\bar{\mathcal{A}}}$. Thus in this paper, we invest methods which can also be applied into the non-trivial case when $0<\kappa< m-\sigma$. One way to achieve a resilient convex combination is through Tverberg points as in \cite{mendes2015DC,LN14ACM,N14ICDCN}.
For any $\mathcal{S}\subset \mathcal{A}$, let $\mathcal{H}(x_\mathcal{S})$ denote the convex hull of vectors in $x_\mathcal{S}$, that is,
\begin{align}
\label{HSdef}\mathcal{H}(x_\mathcal{S})=\{\sum_{k=1}^{|x_\mathcal{S}|}\alpha_k s_k:s_k\in x_\mathcal{S}, \alpha_k\geq 0, \sum_{k=1}^{|x_\mathcal{S}|}\alpha_k=1\}.
\end{align}
Then the existence of Tverberg points is guaranteed by the following theorem: \\

\noindent\textbf{\textit{Tverberg Theorem}} \cite{Tverberg66}:
Suppose $m\ge\kappa(n+1)+1$ for the given set $x_{\mathcal{A}}$.  Then there must exist a partition of $\mathcal{A}$ into $\kappa+1$ disjoint subsets ${\mathcal{B}_1},\cdots,{\mathcal{B}_{\kappa+1}}$ such that
\begin{align}\label{deftver}
\mathcal{T}=\bigcap_{j=1}^{\kappa+1}\mathcal{H}(x_{\mathcal{B}_j})\neq \emptyset,
\end{align}
where
$$\bigcup_{j=1}^{\kappa+1}\mathcal{B}_j=\mathcal{B} $$ and
$$\quad \mathcal{B}_j \cap\mathcal{B}_k=\emptyset,\quad \forall j\neq k. $$ Points in the non-empty intersection $\mathcal{T}$ in \eqref{deftver} are called \textbf{Tverberg points} of the $(\kappa+1)$-partition of $\mathcal{A}$.
\\

While results in \cite{mendes2015DC,LN14ACM,N14ICDCN,HS17TOR} are elegant, one major concern of applying Tverberg points lies in the requirement of high computational complexity. As mentioned in \cite{WD13DCG,PME08ACM}, except for some specific values of $n$, the computational complexity of calculating Tverberg points grows exponentially with the dimension $n$. In the following, we will develop a low-complexity algorithm for achieving resilient convex combinations based on the intersection of convex hulls.

\section{Resilient Convex Combination }\label{Sec_RCC}

\subsection{ A Resilient Convex Combination through Intersection of Convex Hulls}
Let
\begin{align}\label{wopi}
\mathcal{R}=\bigcap_{j=1}^r\mathcal{H}(x_{\mathcal{A}_j}),
\end{align} where $r=\binom {m-\sigma} {m-\sigma-\kappa}$  and $\mathcal{A}_j$, $j=1,2,...,r$, denote all subsets of $\mathcal{A}$ such that
\begin{align}\label{eq_Mit}
\bar{\mathcal{A}}\subset \mathcal{A}_j\subset \mathcal{A}, \  |\mathcal{A}_j|=m-\kappa.
\end{align}
Then one has the following lemma:
\begin{lemma}\label{Lemma_Resi}
If $\mathcal{R}\neq \emptyset$, then any point in $\mathcal{R}$ is a resilient convex combination.
\end{lemma}
\noindent{\bf Proof of Lemma \ref{Lemma_Resi}:} Since the number of malicious points in $x_\mathcal{A}$ is upper bounded by $\kappa$, there must exist at least one subset $\mathcal{A}_{j^*}$ which consists of only normal points. As long as $\mathcal{R}\neq \emptyset$,  for any vector $q\in\mathcal{R}$, it must be true that $q\in \mathcal{H}(x_{\mathcal{A}_{j^*}})$. Thus $q$ is a resilient convex combination. $\qed$

Compared with the Tverberg points set in (\ref{deftver}), the $\mathcal{R}$ in (\ref{wopi}) defines a larger set for choosing resilient convex combinations, as indicated by the following lemma.
\begin{lemma}\label{Lemma_tver} The set $\mathcal{T}$ in \eqref{deftver} and the set $ \mathcal{R}$ in \eqref{wopi} satisfy
	\begin{align}\label{eq_subsub}
	\mathcal{T}\subset \mathcal{R}.
	\end{align}

\end{lemma}

\noindent{\bf Proof of Lemma \ref{Lemma_tver}:}  We first claim that for each $\mathcal{A}_j$, $j=1,2,...,r$, defined in (\ref{eq_Mit}), one of ${\mathcal{B}_1},\cdots,\mathcal{B}_{\kappa+1}$ must be its subset. We prove this by contradiction. Suppose there exists a $\mathcal{A}_{j^\dagger}$ such that none of $\mathcal{B}_1,\cdots,\mathcal{B}_{\kappa+1}$ is a subset of $\mathcal{A}_{j^\dagger}$. Then each $\mathcal{B}_j$, $j=1,2,\cdots,\kappa+1$, must have at least one element that is not in $\mathcal{A}_{j^\dagger}$. Note that any two of ${\mathcal{B}_1},\cdots,{\mathcal{B}_{\kappa+1}}$ are disjoint. Then there are at least $\kappa+1$ elements that are not in $\mathcal{A}_{j^\dagger}$. Then $|\mathcal{A}_{j^\dagger}|\leq m-\kappa-1$, which contradicts the fact that $|\mathcal{A}_{j^\dagger}|=m-\kappa$. Thus for each $\mathcal{A}_j$, one of ${\mathcal{B}_1},\cdots,{\mathcal{B}_{\kappa+1}}$ must be a subset of $\mathcal{A}_j$. From this and the definition of $\mathcal{R}$ in (\ref{wopi}), one has $\bigcap_{j=1}^{\kappa+1}\mathcal{H}(x_{\mathcal{B}_j})\subset \mathcal{R}$, which is (\ref{eq_subsub}).
We complete the proof. \qed
\smallskip

To guarantee that $\mathcal{R}$ is not empty, one has the following lemma:
\begin{lemma}\label{Lemma_nonem}
If $\bar{\mathcal{A}}\neq \emptyset$, then $\mathcal{R}\neq \emptyset$; If $\bar{\mathcal{A}}=\emptyset$, but $m\ge(\kappa(n+1)+1)$, then $\mathcal{R}\neq \emptyset$.
\end{lemma}

\noindent{\bf Proof of Lemma \ref{Lemma_nonem}:}
For the case of $\bar{\mathcal{A}}\neq \emptyset$, recall equation \eqref{eq_Mit} that $\bar{\mathcal{A}}\subset \mathcal{A}_j$. Then for any $j=1,2,...,r$, one has
$$\mathcal{H}(x_{\bar{\mathcal{A}}})\subset\mathcal{H}(x_{\mathcal{A}_j}).$$
It follows that
$$\mathcal{H}(x_{\bar{\mathcal{A}}})\subset\bigcap_{j=1}^r\mathcal{H}(x_{\mathcal{A}_j})=\mathcal{R}$$
Thus, $\bar{\mathcal{A}}\neq \emptyset$ leads to  $\mathcal{R}\neq \emptyset$.

For the case of $\bar{\mathcal{A}}= \emptyset$. If $m\ge(\kappa(n+1)+1)$, one has from Tverberg Theorem that $\mathcal{T}\neq\emptyset$. Thus, Lemma \ref{Lemma_tver} leads to
$$\mathcal{T}\subset\mathcal{R}\neq \emptyset$$
This completes the proof. \qed


\subsection{A Low-Complexity Algorithm to Calculate $\mathcal{R}$}
Since any point in $\mathcal{R}$ is a resilient convex combination (by Lemma \ref{Lemma_Resi}), it is desirable to calculate the set $\mathcal{R}$, which by (\ref{wopi}) is the intersection of a group of convex hulls. Existing approaches for the computation of intersection of convex hulls are usually computationally complex (\#p-hard in \cite{L07PCG,CDH96ConvexH}). Thus in this section we will develop an algorithm with low computational complexity for calculating a point in $\mathcal{R}$.

First, we will propose an equivalent expression of the set $\mathcal{R}$ in terms of equality and inequality \textbf{constraints}.
For each $\mathcal{A}_j=\{j_1,j_2,\cdots,j_p\}$, $j=1,2,...,r$, we define the following matrix 	
\begin{align}\label{eq_y0}
Y_j=\begin{bmatrix}
x_{j_1}(t)&x_{j_2}(t)&\cdots&x_{j_p}(t)
\end{bmatrix}\in\mathbb{R}^{n\times p}. 
\end{align}
We call
\begin{align}
\label{eq_XX}X={\rm diag}\{Y_j,~j=1,2,...,r\}\in\mathbb{R}^{nr\times pr}
\end{align}
the \emph{coordinate matrix}. For example, suppose $\mathcal{A}=\{1,~2,~3\}$, $\bar{\mathcal{A}}=\{1\}$ and $\kappa=1$, then $p=2$, $r=2$ and one has:
\begin{align*}
\mathcal{A}_1&=\{1,2\},\quad \mathcal{A}_2=\{1,3\}\\
Y_1&=\begin{bmatrix}
x_1&x_2
\end{bmatrix}, \quad
Y_2=\begin{bmatrix}
x_1&x_3
\end{bmatrix}\\
X&=\begin{bmatrix}
Y_1&0\\
0& Y_2
\end{bmatrix}.
\end{align*}
This coordinate matrix allows us to characterize the set $\mathcal{R}$ with the following lemma
\begin{lemma} \label{TQP}
	Let $C\in \mathbb{R}^{r\times r}$ be the circulant matrix with the first row in the form of
	$\begin{bmatrix}1&-1&0&\cdots &0\end{bmatrix}$.	
	Then
	\begin{align}\label{Const_ConvH}
	{\mathcal{R}}=\left\{ \frac{1}{r}(\bm{1}'_{r}\otimes I_n)X\bm{\beta}\right\},
	\end{align}	
	for all $\bm{\beta}\in\mathbb{R}^{pr}$ that satisfies
	\begin{align}
	(C\otimes I_n) X\bm{\beta}&=0 \label{eq_b1}\\
	(I_{r}\otimes \bm{1}'_p)\bm{\beta}&=\bm{1}_{r} \label{eq_b2} \\
	\bm{\beta}&\ge 0.\label{eq_b3}
	\end{align}
	
\end{lemma}

\smallskip

Before proving Lemma \ref{TQP}, note that if we let $\bm{\beta}=\col\{\beta_j,~j=1,2,...,r\}$ and $y_j=Y_j\beta_j\in\mathbb{R}^n$, with $\beta_j\in \mathbb{R}^p $, then
\begin{align}\label{yeq}
X\bm{\beta}=\col\{y_j,~j=1,2,...,r\}.
\end{align}
This indicates that the components $y_j$ of $X\bm{\beta}$ are linear combinations of vectors stored in $Y_j$ according to the coefficient vector $\beta_j$.
With this, we carry out the following proof.

\medskip
\noindent{\bf Proof of Lemma \ref{TQP}:}
From \eqref{yeq}, the definition of the circulant matrix $C$, and $(C\otimes I_n) X\bm{\beta}=0$ in (\ref{eq_b1}), one has for all $j=1,2,...,r$, there exists a $y^*$ such that
\begin{align}\label{eq_equal0}
y_j=y^*.
\end{align}

Let $\mathcal{Y}_j$ denote the set of $y_j$ for all $\bm{\beta}$ satisfying (\ref{eq_b2})-(\ref{eq_b3}). Note that these equations ensure $\beta_j$ is a nonnegative vector with entries summing to 1, which guarantees the combinations $y_j=Y_j\beta_j$ are convex. Then, given the definitions \eqref{HSdef} and \eqref{eq_y0}, it is true that
\begin{align}\label{eq_Yk}
\mathcal{Y}_j=\mathcal{H}(x_{\mathcal{A}_j}),\quad j=1,2,...,r.
\end{align}
This along with \eqref{eq_equal0} indicate that the set of all feasible $y^*$ is given by
\begin{align}
\{y^*\}=\bigcap_{j=1}^{r}\mathcal{Y}_j=\bigcap_{j=1}^{r}\mathcal{H}(x_{\mathcal{A}_j})=\mathcal{R}
\end{align}
Recall that $(\bm{1}'_{r}\otimes I_n)X\bm{\beta}=y_1+y_2+...+y_r=ry^*$. Thus,
$$\left\{ \frac{1}{r}(\bm{1}'_{r}\otimes I_n)X\bm{\beta}\right\}=\mathcal{R}.$$
This completes the proof.  \qed

\medskip

Lemma \ref{TQP} tells us that computing the set $\mathcal{R}$ is equivalent to solving equations \eqref{Const_ConvH}-\eqref{eq_b3}. However, to obtain a particular resilient convex combination, one does not necessarily have to find all points in this set. Here, we consider a point $u$, which tends to equally use all vectors in $x_{\mathcal{A}}$.
\begin{align} \label{eq_bvv}
u=\frac{1}{r}(\bm{1}'_{r}\otimes I_n)X\bm{\beta}^*
\end{align}
where $\bm{\beta}^*$ is computed by the following quadratic programming problem:
\begin{align}
\text{minimize} &\quad J(\bm{\beta})=\frac{1}{p}\|\bm{\beta}-\frac{1}{p}\bm{1}_{pr}\|_2^2 \label{eq_opti}
\end{align}
subject to  constraints \eqref{eq_b1}-\eqref{eq_b3}. Specially when $\kappa=0$,  one has $r=1$, $p=m$. Then $\bm{\beta^*}=\frac{1}{m}\bm{1}_{m}$, and
$u=\frac{1}{m}\sum_{j=1}^m x_j$, which is the average of all vectors in $x_\mathcal{A}$.

\medskip

To reveal the mechanism of \eqref{eq_opti}, note that
\begin{align*}
\frac{1}{p}\|\bm{\beta}-\frac{1}{p}\bm{1}_{pr}\|_2^2=\frac{1}{p}\sum_{j=1}^r\|\beta_j-\frac{1}{p}\bm{1}_{r}\|_2^2
\end{align*}
which minimizes the sum of variances of coefficient vectors $\beta_j$.  Recall that the entries of $\beta_{j}$ summing up to 1, its average should be $\frac{1}{p}$ so when the value of $J(\bm{\beta})$ approaches $0$, the weights are equally distributed to all neighbors' states. In this way, the $u$ in \eqref{eq_bvv} can be viewed as an unbiased choice of the resilient convex combination that lies in the region $\mathcal{R}$.
It is worth mentioning that the complexity of achieving the $\bm\beta^*$ in (\ref{eq_opti}) is $\mathcal{O}(n\cdot (mr)^3)$\cite{SL04ConvexOpt}.


\subsection{Main Result and Comparison}\label{Sec_Tve} The main result of the paper is the following theorem, which summarizes Lemmas \ref{Lemma_Resi}-\ref{TQP}:
\begin{theorem}\label{Th1}
Consider a set $x_{\mathcal{A}}=\{x_1,x_2,\cdots,x_m\}$ of $m$ vectors in $\mathbb{R}^n$, where $\mathcal{A}=\{1,2,...,m\}$. Suppose at most $\kappa$ vectors in $x_{\mathcal{A}}$ are malicious and $\bar{\mathcal{A}}$ is a known label set that only contains normal vectors, which can also be empty. If $\bar{\mathcal{A}}\neq\emptyset$ or $m\ge (\kappa(n+1)+1)$, any point in the non-empty set $\mathcal{R}$ defined in \eqref{wopi} is a resilient convex combination. One specific point $u\in \mathcal{R}$ defined in \eqref{eq_bvv} can be computed by solving the quadratic programming problem \eqref{eq_opti} subject to  constraints \eqref{eq_b1}-\eqref{eq_b3}.
\end{theorem}

Recall that the set of Tverberg points $\mathcal{T}$, which is a subset of $\mathcal{R}$, also provides a region for choosing resilient convex combinations. Comparisons between $\mathcal{T}$ and $\mathcal{R}$ are provided as follows:
\begin{itemize}
	\item First, determining Tverberg points requires \textbf{high computational complexity}. Although the Tverberg Theorem provides a sufficient condition for the existence of a partition leading to Tverberg point, the theorem does not provide an algorithm for finding the partition for achieving Tverberg points, apart from enumerating all possible partitions and checking the intersection of convex hulls for each partition. As mentioned in \cite{WD13DCG,PME08ACM}, except for some specific values of $n$, the computational complexity of achieving Tverberg points grows exponentially with the dimension $n$.
	In contrast, it has been shown that the resilient convex combination proposed in (\ref{eq_bvv}) can be computed by solving a standard quadratic programming problem, whose computational complexity is polynomial in $n$.
	
	\item Second, the existence of an non-empty Tverberg point set $\mathcal{T}$ requires that $m\geq \kappa(n+1)+1$, while one has $\mathcal{R}\neq\emptyset$ if $\bar{\mathcal{A}}\neq\emptyset$ or $m\geq \kappa(n+1)+1$. In achieving resilience for distributed algorithms, one aims to guarantee all normal agents' states to converge to a consensus. Then each normal agent at least has one element (which is itself) in $\bar{\mathcal{A}}$. Then the existence of an non-empty $\mathcal{R}$ is automatically guaranteed. Please refer to Fig. \ref{IntConv}. (A) and (B) for an example when $\mathcal{R}$ is non-empty while $\mathcal{T}=\emptyset$, and an example $\mathcal{T}\subset \mathcal{R}$, respectively.
	
	\begin{figure}[h]
		\centering
		\subfloat[$\kappa=2$, $m=6$, $ m<\left(\kappa(n+1)+1\right)$]{{\includegraphics[width=5cm]{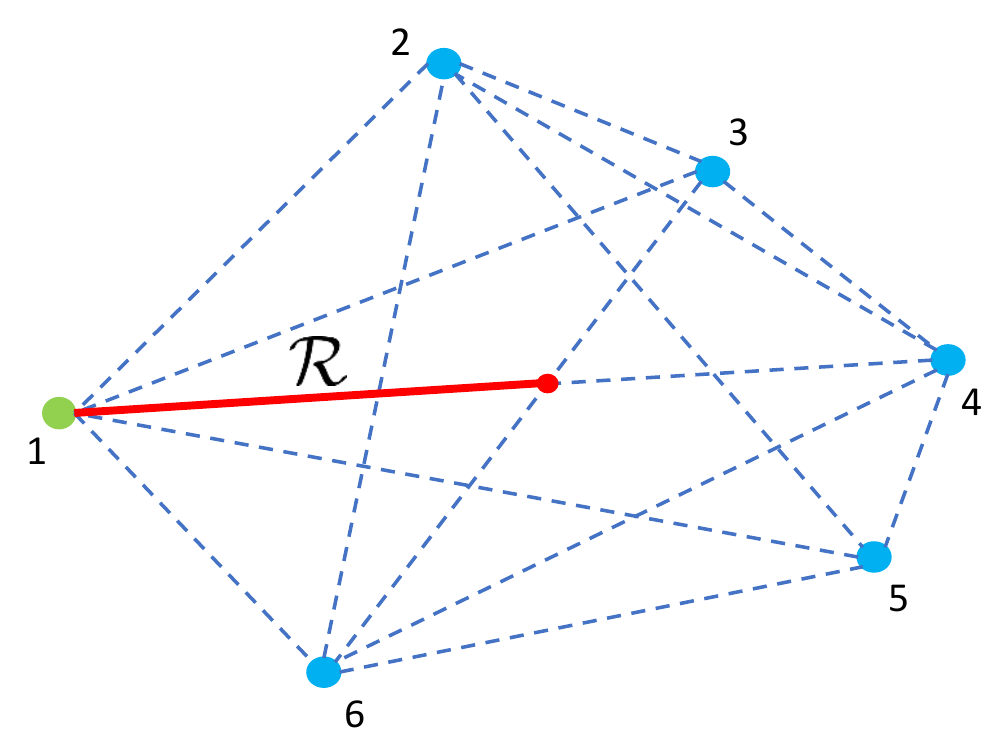} }}%
		\qquad
		\subfloat[$\kappa=1$, $m=4$, $m=\left(\kappa(n+1)+1\right)$]{{\includegraphics[width=5cm]{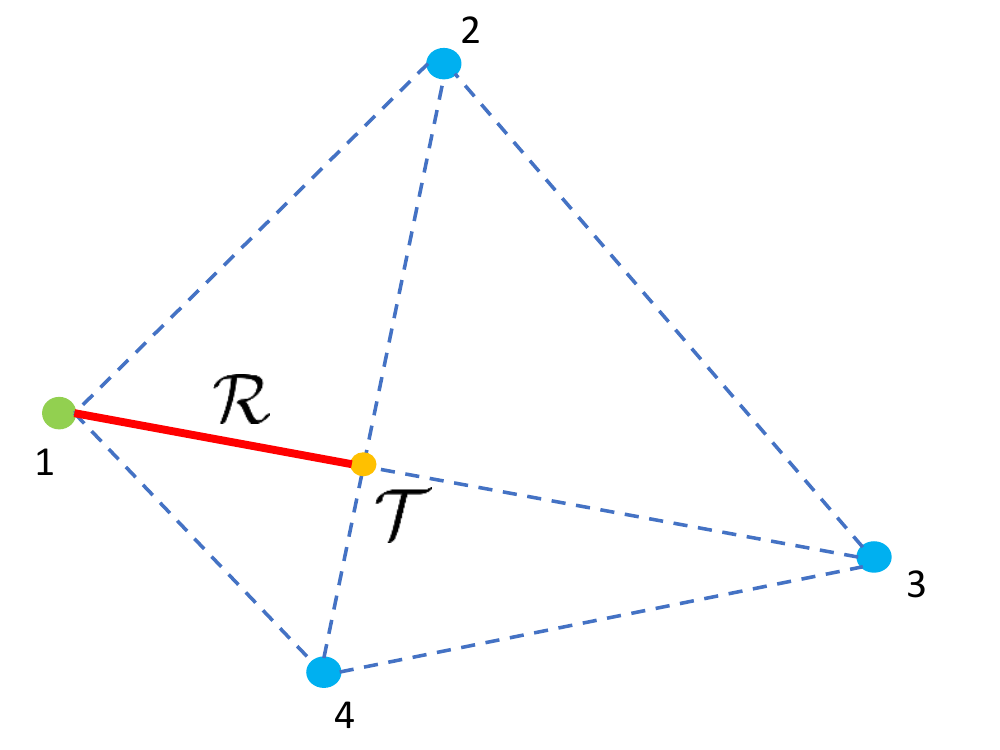} }}%
		\caption{Finding Tverberg point $\mathcal{T}$ (yellow)  and $\mathcal{R}$ (red) in a 2-D space, with $\bar{\mathcal{A}}=\{1\}$.}%
		\label{IntConv}%
		
		
	\end{figure}
	
	
\end{itemize}


\section{Application of the Resilient Convex Combination into Consensus-Based Distributed Algorithms}\label{Sec_App}
%
Consider a network of $\bar{m}$ agents in which each agent $i$ is able to sense or receive information from certain other nearby agents, termed agent $i$'s neighbors. We suppose agent $i$ is always a neighbor of itself and we let $\mathcal{N}_i(t)$ denote the set of agent $i$'s neighbors at time $t$, $i=1,2,\cdots,\bar{m}$. The neighbor relations can be described by a time-dependent graph $\mathbb{G}(t)$ such that there is a directed edge from $j$ to $i$ in $\mathbb{G}(t)$ if and only if $j\in \mathcal{N}_i(t)$. Suppose each agent $i$ controls a state vector $x_i(t)\in \mathbb{R}^n$. Consensus-based distributed algorithm solves consensus problems that are unconstrained\cite{MAB08SIAM}, constrained by linear or nonlinear constraints\cite{SJA15TAC}, and/or minimize a global objective function\cite{AAP10TAC}. These algorithms share a common form
\begin{align}
x_i(t+1)
=f_i(x_i(t), v_i(t)) \label{mainalgorithm}
\end{align}
	 where $v_i(t)$ is a convex combination of all agent $i$'s neighbors' states, that is,
$$v_i(t)=\sum\limits_{j\in\mathcal{N}_i(t)} w_{ij}(t)x_j(t)$$ with $\sum\limits_{j\in\mathcal{N}_i} w_{ij}(t)=1$.

Since each agent updates its state by a convex combination of all its neighbors' states, when one or more neighbors are malicious, the convex combination also contains false information, which may prevent the overall consensus goal from being reached. This motivates the \textbf{key idea} to replace the convex combination $v_i(t)$ with a \textbf{resilient convex combination} $u_i(t)$ defined in \eqref{eq_bvv}. To be more specific, in each time step, we assume that agent $i$ knows $\mathcal{A}=\mathcal{N}_i(t)$, $\bar{\mathcal{A}}=\{i\}$ and the upper bound of agent's malicious neighbors $\kappa$. Then it computes $u_i(t)$ by solving the quadratic programming problem \eqref{eq_opti} under constraints \eqref{eq_b1}-\eqref{eq_b3}. In this way, the malicious information is \textbf{automatically} isolated by $u_i(t)$.

\subsection{Simulation}\label{Sec_Simu}

\begin{figure}[h]
	\centering
	\includegraphics[width=8 cm]{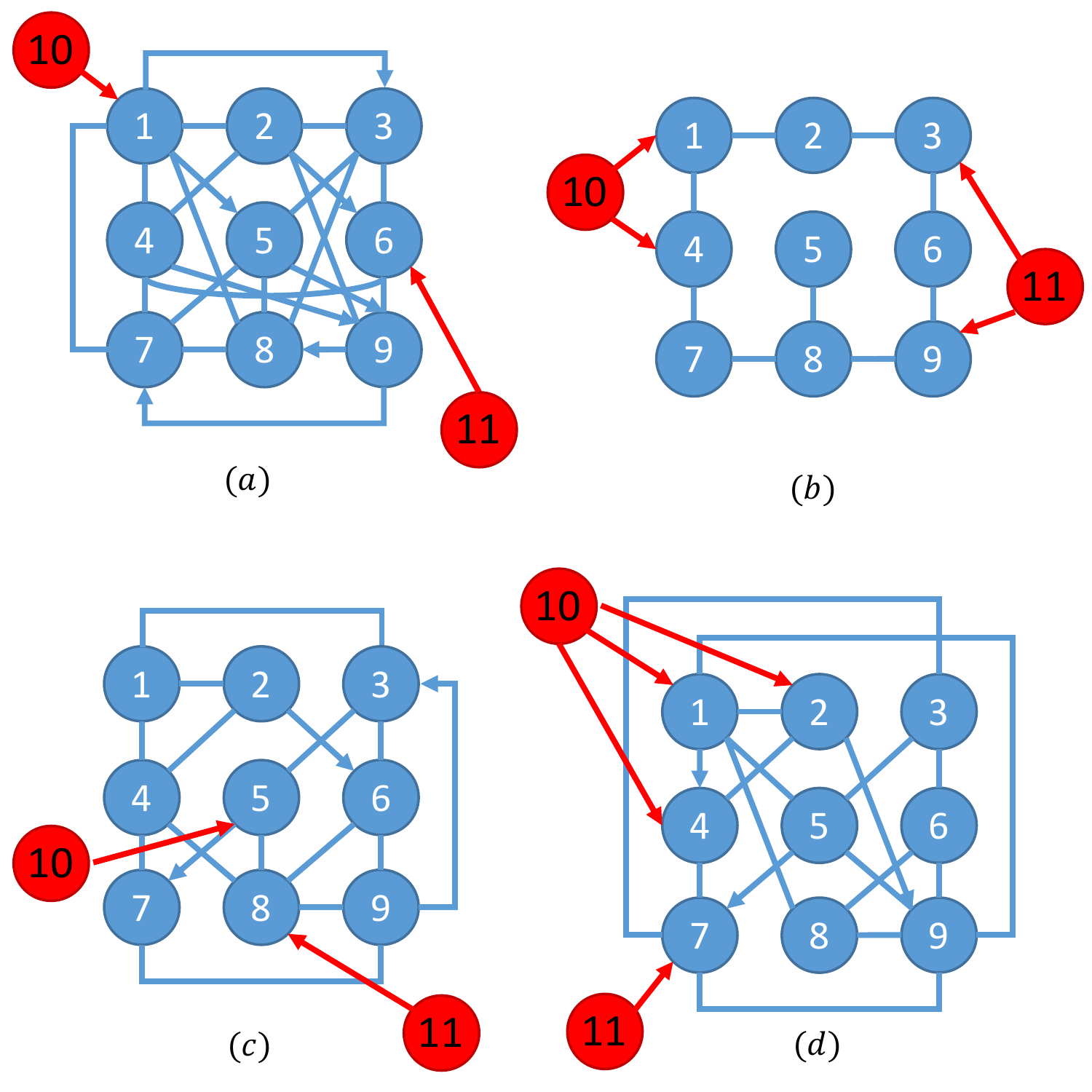}
	\caption{A network of 11 agents with malicious agents marked in red.}
	\label{Fig1}
\end{figure}

In this section, we provide simulations for a 11-agent time-varying network consisting of both directed and undirected edges as indicated in Fig. \ref{Fig1}, in which agent $10$ and $11$ are malicious agents and connect themselves to different normal agents as time evolves. By replacing $v_i(t)$ in (\ref{mainalgorithm}) with $u_i(t)$, the problem of interest is to check whether all normal agents from $1$ to $9$ under this update  still reach the desired consensus in the presence of malicious agents.  In the following examples, we suppose each $x_i(t)\in \mathbb{R}^2$, $\kappa=1$ (even though there are two malicious agents, for each agent, the upper bound of malicious neighbor is 1). Each malicious agent sends a state to its neighbors which is randomly chosen from the set $[0,2]\times[0,2]$.

\noindent{\bf Example 1 (Unconstrained Consensus).}

We first consider the unconstrained consensus problem, in which all $x_i(t)\in \mathbb{R}^2$, $i=1,2,\cdots,9$, aim to reach consensus by the following update
\eq{x_i(t+1)=\frac{1}{d_i(t)}\sum\limits_{j\in\mathcal{N}_i(t)}x_i(t)\label{eq_consensus}} with initialization $x_i(0)$ randomly chosen in from the areas of $[0,2]\times[0,2]$.
Let
\begin{align}\label{defV}
V(t)=\frac{1}{2}\sum\limits_{i=1}^8\|x_i(t)-x_{i+1}(t)\|_2^2
\end{align} which measures the closeness of all normal agents' states to consensus.

\smallskip

\noindent \textit{A. Under Fixed Graph}

\begin{figure}[h]
	\centering
	\includegraphics[width=8 cm]{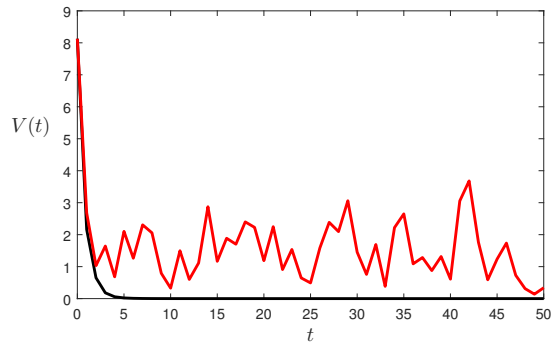}
	\caption{Simulations of normal agents under the consensus update (\ref{eq_consensus}) without malicious agents (blank line) and with malicious agents $10$ and $11$ (red line).}
	\label{Fig2}
\end{figure}

Suppose the network is a fixed one as in Fig. \ref{Fig1}-(a).
Simulation results for the consensus update (\ref{eq_consensus}) with and without the presence of malicious agents are shown in Fig. \ref{Fig2}, which indicates unsurprisingly that the traditional consensus update (\ref{eq_consensus}) could easily fail in the presence of malicious agents.

\begin{figure}[h]
	\centering
	\includegraphics[width=8 cm]{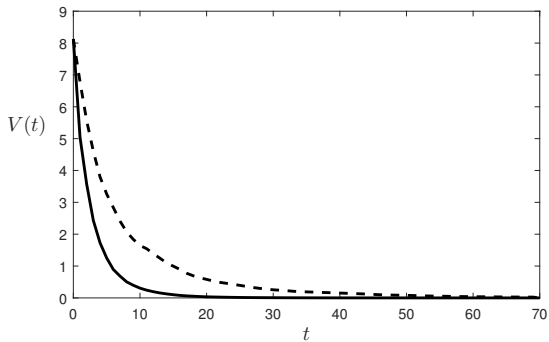}
	\caption{Consensus is reached by introducing $u_i(t)$ as Tverberg points (indicated by the dash line) or as the resilient convex combination  (\ref{eq_bvv}) (indicated by the solid line).}
	\label{Fig3}
\end{figure}

By introducing a new resilient convex combination ${u}_i(t)$ at each agent $i$, which is a convex combination of normal agents, one could employ the following update
\eq{\label{eq_rcon}x_i(t+1)=u_i(t)}
where $u_i(t)$ could be chosen as Tverberg points or as the resilient convex combination in (\ref{eq_bvv}). Consensus is reached in the presence of malicious agents in both cases as shown by simulations in Fig. \ref{Fig3}. It is also worth mentioning that using $u_i$ defined in (\ref{eq_bvv}) to replace the original convex combination $v_i(t)$ could lead to faster convergence than using Tverberg points as also indicated in Fig. \ref{Fig3}
\smallskip

\noindent \textit{B. Under Time-varying Graph}

\begin{figure}[h]
	\centering
	\includegraphics[width=8 cm]{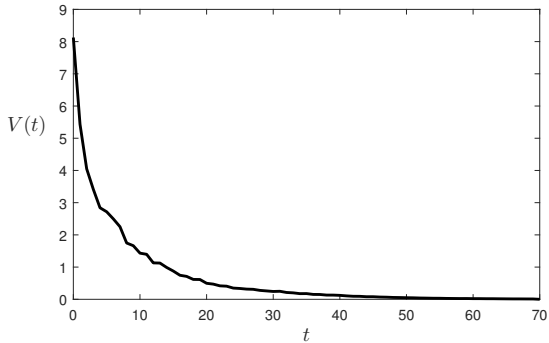}
	\caption{Simulations by using the resilient convex combination $u_i(t)$ of (\ref{eq_bvv}) into (\ref{eq_rcon}).}
	\label{Consensus2TV}
\end{figure}

We perform simulations of unconstrained consensus on a periodic sequence of time-varying networks as in Fig. \ref{Fig1}.  The method based on Tverberg point is not applicable here since the number of each agent's neighbors is not always greater than the condition required by the Tverberg Theorem. However, one could still reach consensus by introducing the resilient convex combination  (\ref{eq_bvv}) into (\ref{eq_rcon}). Simulation results are shown in Fig. \ref{Consensus2TV}.

\medskip
\noindent{\bf Example 2 (Constrained Consensus). }

\begin{figure}[h]
	\centering
	\includegraphics[width=8 cm]{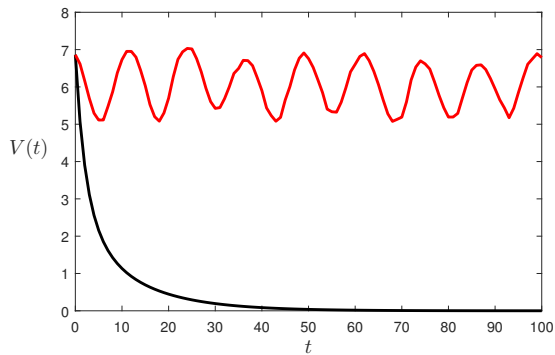}
	\caption{Simulation results under the update (\ref{eq_LE}) with no malicious agents (indicated by the black line) or with malicious agents (indicated by the red line). }
	\label{Fig4}
\end{figure}

\begin{figure}[h]
	\centering
	\includegraphics[width=8 cm]{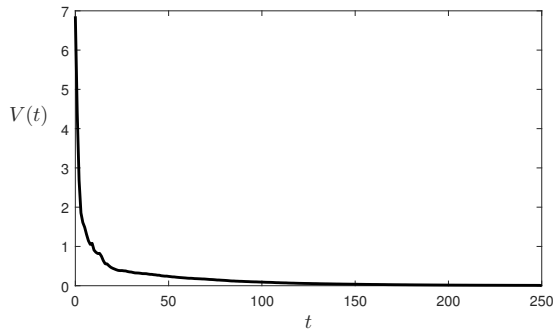}
	\caption{Simulations by using the resilient convex combination $u_i(t)$ of (\ref{eq_bvv}) in (\ref{eq_LE}). }
	\label{Fig5}
\end{figure}

We consider the distributed algorithm for solving linear equations. Suppose each agent $i$ knows
\begin{align*}
A_ix_i=b_i
\end{align*}
where
\begin{align*}
&A_1=A_2=A_3=\begin{bmatrix} 3&-1 \end{bmatrix} &b_1=b_2=b_3=2\\
&A_4=A_5=A_6=\begin{bmatrix} 0& 1 \end{bmatrix} &b_4=b_5=b_6=1\\
&A_7=A_8=A_9=\begin{bmatrix} -1&3 \end{bmatrix} &b_7=b_8=b_9=2
\end{align*}
and updates its state according to
\eq{x_i(t+1)=x_i(t)-P_i(x_i(t)-v_i(t))\label{eq_LE}}
where $P_i$ is the orthogonal projection on the kernel of $A_i$ and $v_i(t)=\frac{1}{d_i(t)}\sum\limits_{j\in\mathcal{N}_i}x_i(t)$. Simulations are still performed on a periodic sequence of time-varying networks as in Fig. \ref{Fig1}. Let $x^*=\begin{bmatrix}
1&1&1
\end{bmatrix}'$ denote a solution to $Ax=b$, which is also the desired consensus value. Let $$
V(t)=\frac{1}{2}\sum\limits_{i=1}^9\|x_i(t)-x^*\|_2^2
$$ which measures the closeness between all agents' states to $x^*$. Simulation results are as shown in Fig. \ref{Fig4} for the cases with and without malicious agents, respectively. The presence of malicious agents also disrupts the distributed algorithm (\ref{eq_LE}) for solving linear equations.
By using the resilient convex combination $u_i(t)$ of (\ref{eq_bvv}) at each agent $i$ in (\ref{eq_LE}), one still enables all agents to achieve $x^*$ exponentially fast in the presence of malicious agents as shown in Fig. \ref{Fig5}.

\section{Conclusion}\label{Conc}
Given a set of vectors that includes both normal and malicious information, this paper has proposed a way to determine a resilient convex combination
by using intersection of convex hulls. By formulating the set of such combinations as linear constraints, a vector inside this set can be computed by quadratic programming.  It has been shown that the obtained resilient convex combination  can  isolate harmful state information injected by cyber-attacks.
In addition, since no identification process is required, the method has promise for dealing with time-varying attacks for consensus-based distributed algorithms, as shown by simulations.

\bibliographystyle{AIMS}

\bibliography{ConsensusGossip,CyberSecurity,FormationControl,OptiCompu,Others,Shaoshuai,Anderson}

\medskip

\medskip

\end{document}